\documentclass{amsart}

\usepackage{amssymb}
\usepackage{amsmath}
\usepackage{amsthm}
\usepackage[mathscr]{eucal}

\newtheorem{thrm}{Theorem}[section]
\newtheorem{lemma}[thrm]{Lemma}
\newtheorem{prop}[thrm]{Proposition}

\theoremstyle{definition}

\theoremstyle{remark}
\newtheorem{remark}[thrm]{Remark}

\numberwithin{equation}{section}

\begin{document}

\bibliographystyle{plain}

\title{Solution of the $\bar{\partial}$-Neumann problem on a bi-disc}

\author{Dariush Ehsani}

\address{Department of Mathematics, Texas A\&M University, College Station,
Texas 77843-3368} \email{ehsani@math.tamu.edu}


\begin{abstract}
In this paper we study the behavior of the solution to the
$\bar{\partial}$-Neumann problem for $(0,1)$-forms on a bi-disc in
$\mathbb{C}^2$.  We show singularities which arise at the
distinguished boundary are of logarithmic and arctangent type.
\end{abstract}

\maketitle

\section{Introduction}
Let $\Omega\subset\mathbb{C}^n$ be a bounded, pseudoconvex domain,
equipped with the standard Hermitian metric. The
$\bar{\partial}$-Neumann problem, on domains with a $C^2$ defining
function, takes the form of the boundary value problem
\begin{equation*}
\square u=f \mbox{ in } \Omega \nonumber,
\end{equation*}
for $f$ in $L^2_{p,q}(\Omega)$, and
\begin{align}
u \rfloor \bar{\partial}\rho &= 0,  \\
\bar{\partial}u \rfloor \bar{\partial}\rho &= 0 \label{bndry},
\end{align}
on $\partial \Omega$, where $\square$ is the complex Laplacian,
$\bar{\partial} \bar{\partial}^{\ast} + \bar{\partial}
\bar{\partial}^{\ast}$.

In the past decade, considerable attention has been given to the
study of the $\bar{\partial}$-Neumann problem on non-smooth
domains.  We point to the papers of Henkin and Iordan \cite{HI},
Henkin, Iordan, and Kohn \cite{HIK}, Michel and Shaw
\cite{MS1,MS2}, and Straube \cite{St}, in which properties,
compactness and subelliptic estimates, hold for the Neumann
operator, N, the inverse to the $\bar{\partial}$-Neumann problem,
on certain non-smooth domains.

In \cite{Eh}, the author studied the $\bar{\partial}$-Neumann
problem for $(0,1)$-forms on a model domain, the product of two
half-planes in $\mathbb{C}^2$.  We continue here the study of the
problem for $(0,1)$-forms on model domains, focusing on the
bi-disc, $\Omega =\mathbb{D}_1\times\mathbb{D}_2 \in
\mathbb{C}^2,$ where $\mathbb{D}_1 \subset \mathbb{C}$ and
$\mathbb{D}_2\subset \mathbb{C}$ are defined by the equations $r_1
< 1$ and $r_2 < 1$, respectively, where $r_j=|z_j|$, $j=1,2.$ The
existence of a solution in $L^2(\Omega)$ is given by H\"{o}rmander
\cite{Hor1}. We shall see singularities only occur on the
distinguished boundary, $\partial
\mathbb{D}_1\times\partial\mathbb{D}_2$. Our main result is the
\begin{thrm}
\label{main} Let $\Omega \in \mathbb{C}^2$ be the bi-disc,
$\mathbb{D}_1\times
 \mathbb{D}_2$, where
$\mathbb{D}_j$ is the disc $\{z_j:|z_j|<1\}$ for $j=1,2$.
  Let $f=f_1d\bar{z}_1
+f_2 d\bar{z}_2$ be a $(0,1)$-form such that $f \in
C^{\infty}_{(0,1)}(\overline{\Omega})$, the family of
$(0,1)$-forms whose coefficients are in
$C^{\infty}(\overline{\Omega})$, and $u=u_1 d\bar{z}_1 +u_2
d\bar{z}_2$ the $(0,1)$-form which solves the
$\bar{\partial}$-Neumann problem with data the $(0,1)$-form $f$ on
$\Omega$. Then, with $z_j=r_je^{i\theta_j}$, near $r_1=r_2=1$,
$u_j$ can be written as
\[
u_j =
 \alpha_j \log \big((\log r_1)^2+(\log r_2)^2\big) + \beta_j +\gamma_j
\arctan \left( \frac{\log r_1}{\log r_2} \right) \qquad j=1,2,
\]
where $\alpha_j$, $\beta_j$, $\gamma_j$ are smooth functions of
$r_1,r_2,\theta_1,\theta_2$.
\end{thrm}

We point out the interesting phenomenon that, although the Neumann
operator on $(0,1)$-forms, $N_{(0,1)}$, is not regular,
$\bar{\partial}^{\ast} N_{(0,1)}$ is.  Regularity of
$\bar{\partial}^{\ast} N_{(0,1)}$ follows from regularity of the
Bergman projection on the bi-disc and the formula for the Bergman
projection, P,
\begin{equation*}
Pg=g-\bar{\partial}^{\ast}N_{(0,1)}\bar{\partial}g
\end{equation*}
for $g\in L^2(\Omega)$.

I wish to thank several people who helped me in the preparation of
this article.  I offer my gratitude to David Barrett under whose
guidance and encouragement I first considered the problem on the
bi-disc.  I also wish to thank Harold Boas, Peter Kuchment, and
Emil Straube with whom I could share and discuss various ideas in
the course of my research.

\section{Setup}
We set up the $\bar{\partial}$-Neumann problem for $(0,1)$-forms
on the bi-disc, $\Omega =\mathbb{D}_1\times\mathbb{D}_2 \in
\mathbb{C}^2$ and prove regularity results away from the
distinguished boundary.
\begin{equation*}
\bar {\partial} \bar {\partial}^{\ast}u + \bar {\partial}^{\ast}
\bar{\partial}u=f
\end{equation*}
gives equations for $u_1$ and $u_2$ based on the Laplacian:
\begin{align}
\Delta u_1 &= -2f_1,\label{laplace} \\
\Delta u_2 &= -2f_2, \nonumber
\end{align}
which, in polar coordinates $(r_1,\theta_1)$, $(r_2,\theta_2)$,
are
\begin{align}
 \label{lap1} & \frac{\partial^2 u_j}{\partial r_1^2} +
\frac{1}{r_1} \frac{\partial u_j}{\partial r_1} + \frac{1}{r_1^2}
\frac{\partial^2 u_j}{\partial
\theta_1^2}+ \\
&  \frac{\partial^2 u_j}{\partial r_2^2} +
\frac{1}{r_2}\frac{\partial u_j}{\partial r_2} + \frac{1}{r_2^2}
\frac{\partial^2 u_j}{\partial \theta_2^2} = -2f_j \nonumber.
\end{align}

The boundary conditions (\ref{bndry}), which were defined for
$C^2$ domains, may be adapted to our case of the bi-disc to yield
the conditions
\begin{align}
    u_1&=0 \qquad \mbox{ when } r_1=1 \label{dirb1}, \\
u_2&=0 \qquad \mbox{ when } r_2=1, \nonumber
\end{align}
and
\begin{equation}
\frac{\partial u_2}{\partial \bar{z}_1}-\frac{\partial
u_1}{\partial \bar{z}_2} =0 \label{neubu}
\end{equation}
when $r_1=1$ or $r_2=1$. However, since $u_1=0$ when $r_1=1$, we
must have $\frac{\partial u_1}{\partial \bar{z}_2} =0$, and on the
boundary $r_1=1$, (\ref{neubu}) is
\begin{equation*}
\frac{\partial u_2} {\partial \bar{z}_1}=0.
\end{equation*}
Similarly, for $r_2=1$, (\ref{neubu}) is
\begin{equation}
\frac{\partial u_1} {\partial \bar{z}_2}=0.  \label{neub1}
\end{equation}

\begin{lemma}
 \label{genbi}
 Let $u$ be a solution to the $\bar{\partial}$-Neumann problem on
$\Omega=\mathbb{D}_1\times\mathbb{D}_2$. Then $u$ is smooth in any
neighborhood, $V\subset\overline{\Omega}$ not intersecting
$\partial \mathbb{D}_1\times\partial\mathbb{D}_2$.
\end{lemma}
\begin{proof}
We consider $u_1$, the solution to equation (\ref{laplace}) with
the boundary conditions given by (\ref{dirb1}) and (\ref{neub1}).

Interior regularity follows from the strong ellipticity of the
Laplacian.

Also, general regularity at the boundary arguments for the
Dirichlet problem can be applied to the case in which $V$ is a
neighborhood such that $V\bigcap\partial \Omega=V\bigcap\partial
\mathbb{D}_1\neq \emptyset$ (see \cite{Fo}).

Lastly, suppose $V$ is a neighborhood such that $V\bigcap\partial
\Omega=V\bigcap\partial \mathbb{D}_2\neq \emptyset$. Define
$v=\frac{\partial u_1}{\partial \bar{z}_2}$ and consider the
related problem
\begin{equation*}
\triangle v = -2\frac{\partial f_1}{\partial \bar{z}_2}
\end{equation*}
on $\Omega$, with the conditions
\begin{eqnarray*}
v=0& & \quad \mbox{on } r_1=0, \\
 v=0& & \quad \mbox{on }
 r_2=0.
\end{eqnarray*}

 We know, from above, that $v$ is smooth on all neighborhoods not
intersecting $\partial\mathbb{D}_1\bigcap\partial\mathbb{D}_2$,
hence in $V$. Let $z'=(z_1',z_2')\in V \bigcap
\partial \mathbb{D}_2$.
We will work in the neighborhood $\mathbb{D}_1\times V_2$, where
$V_2$ is a bounded neighborhood of $z_2'$ in
$\overline{\mathbb{D}}_2$ such that $V_2 \bigcap \mathbb{D}_2$ has
smooth boundary.  Let $\chi \in C^{\infty}_0(\overline{V_2})$ such
that $\chi \equiv 1$ near $z_2'$.
 Define
\begin{equation*}
u'=\frac{1}{2\pi i}\int_{V_2} \frac{\chi(\zeta_2)
v(z_1,\zeta_2)}{\zeta_2-z_2} d\zeta_2 \wedge d\bar{\zeta}_2.
\end{equation*}
 $u'$ has the properties $\frac{\partial u'}{\partial
\bar{z}_2}=v$ near $z'$ and $u'\in C^{\infty}(\mathbb{D}_1\times
\overline{V_2})$ \cite{Be1}.

We define the operators $\triangle_j$ to be
$\frac{\partial^2}{\partial x_j^2}+\frac{\partial^2}{\partial
x_j^2}$ for $j=1,2$. Then, computing $\triangle u'$ in
$\mathbb{D}_1\times V_2$, we find
\begin{equation*}
 \triangle u'=-\frac{1}{\pi i}\int_{V_2} \frac{\chi(\zeta_2)
\frac{\partial f_1}{\partial \bar{\zeta}_2}}{\zeta_2-z_2} d\zeta_2
\wedge d\bar{\zeta}_2 +\phi(z_1,z_2),
\end{equation*}
which is in $C^{\infty}(\mathbb{D}_1\times \overline{V_2})$,
where, with $\rho=|\zeta_2|$,
\begin{align*}
\phi(z_1,z_2)=& \frac{1}{2 \pi i}\int_{V_2} \frac{(\triangle_2'
\chi) v(z_1,\zeta_2)}{\zeta_2-z_2}d\zeta_2\wedge d\bar{\zeta}_2
+\frac{2}{\pi i}\int_{V_2} \frac{\frac{\partial \chi} {\partial
\zeta_2} \frac{\partial v}{\partial \bar{\zeta}_2}+ \frac{\partial
\chi} {\partial \bar{\zeta}_2} \frac{\partial v}{\partial
\zeta_2}}{\zeta_2-z_2}
d\zeta_2\wedge d\bar{\zeta}_2\\
&-\frac{1}{2\pi i}\int_{\partial V_2 \bigcap \partial
\mathbb{D}_2} \frac{\chi(\zeta_2)\frac{\partial v(z_1,\zeta_2)
}{\partial \rho}} {\zeta_2-z_2}d\zeta_2.
\end{align*}

We set $w=u_1-u'$ and show $w\in C^{\infty}(\mathbb{D}_1\times
\overline{V_2}).$ For $z_2$ near $z_2'$, $\frac{\partial
w}{\partial \bar{z}_2}=0$, in which case
\begin{equation}
\label{triw}
 \triangle_1 w = \triangle w = -2f_1 - \triangle u'=-2f_1+\frac{1}{\pi i}
 \int_{V_2}
\frac{\chi(\zeta_2) \frac{\partial f_1}{\partial
\bar{\zeta}_2}(z_1,\zeta_2)}{\zeta_2-z_2} d\zeta_2 \wedge
d\bar{\zeta}_2 -\phi(z_1,z_2).
\end{equation}
We also have the boundary condition $w=0$ when $r_1=0$.  Hence $w$
is the solution to a Dirichlet problem on the unit disc,
\begin{equation}
\label{nicesoln}
 w=\int_{\mathbb{D}_1}G_1(z_1,\zeta_1)
\Phi(\zeta_1,z_2) d\zeta_1\wedge d\bar{\zeta}_1,
\end{equation}
 where $G_1$ is the Green's function for
$\mathbb{D}_1,$
\[
G_1=\frac{1}{2\pi} \log|z_1-\zeta_1| - \frac{1}{2\pi} \log \left|
|z_1|^{-1}z_1-|\zeta_1|^{-1}\zeta_1\right|,
\]
 and
$\Phi$ is defined to be the right hand side of equation
\ref{triw}.  Because $\Phi \in C^{\infty}(\mathbb{D}_1\times
\overline{V_2})$, so is $w$, and $u_1\in
C^{\infty}(\mathbb{D}_1\times \overline{V_2})$ follows from the
fact that $u'\in C^{\infty}(\mathbb{D}_1\times \overline{V_2})$.

The same reasoning applies to $u_2$, and this proves the lemma.
\end{proof}

We may simplify our calculations if we consider the equations
\begin{equation*}
\triangle v_i=g_i
\end{equation*}
with boundary conditions
\begin{equation*}
v_i=0 \qquad \mbox{ on } \partial \Omega
\end{equation*}
for $i=1,2$, where $v_i=\frac{\partial u_i}{\partial \bar{z}_j}$
and $g_i=-2\frac{\partial f_i}{\partial \bar{z}_j}$ ($j\neq i$).

Expanding $v_1$ and $g_1$ into Fourier series:
\begin{align}
 v_1&=\sum_{m_1,m_2=-\infty}^{\infty}
a_{m_1m_2}(r_1,r_2)e^{im_1\theta_1}
e^{im_2\theta_2} \label{vexpn}\\
g_1&=\sum_{m_1,m_2=-\infty}^{\infty} \nonumber
c_{m_1m_2}(r_1,r_2)e^{im_1\theta_1} e^{im_2\theta_2}.
\end{align}
Using these expansions in (\ref{lap1}),(\ref{dirb1}), and
(\ref{neub1}), we see
  the
family of equations
\begin{align*}
& \frac{\partial^2 a_{m_1m_2}}{\partial r_1^2} +\frac{1}{r_1}
\frac{\partial a_{m_1m_2}}{\partial r_1} -
\frac{m_1^2}{r_1^2}a_{m_1m_2}+\nonumber \\
& \frac{\partial ^2 a_{m_1m_2}}{\partial r_2^2}
+\frac{1}{r_2}\frac{\partial a_{m_1m_2}}{\partial
r_2}-\frac{m_2^2}{r_2^2}a_{m_1m_2}=c_{m_1m_2} \qquad
m_1,m_2=0,\pm1,\ldots
\end{align*}
 are satisfied with the boundary conditions
\begin{align*}
a_{m_1m_2}(1,r_2) &= 0, \\
 a_{m_1m_2}(r_1,1)&=0.
\end{align*}
We have analogous equations for $v_2$.
\section{Solution}
 We are then led to study the equations
\begin{multline}
\frac{\partial^2 a_{m_1m_2}}{\partial r_1^2} +\frac{1}{r_1}
\frac{\partial a_{m_1m_2}}{\partial r_1} -
\frac{m_1^2}{r_1^2}a_{m_1m_2} +\\
\frac{\partial ^2
a_{m_1m_2}}{\partial r_2^2} +\frac{1}{r_2}\frac{\partial
a_{m_1m_2}}{\partial r_2}-\frac{m_2^2}{r_2^2}a_{m_1m_2}=c_{m_1m_2}
\label{lap}
\end{multline}
in the space $r_1 < 1$, $r_2< 1$. Here
$a_{m_1m_2}=a_{m_1m_2}(r_1,r_2)$ and
$c_{m_1m_2}=c_{m_1m_2}(r_1,r_2),$ and the boundary conditions
\begin{align*}
 a_{m_1m_2}(1,r_2) &= 0 ,\\
a_{m_1m_2}(r_1,1) &= 0
\end{align*}
hold.

We make the transformation $y_j=-\log r_j$ for $j=1,2$ in
(\ref{lap}), and multiply the resulting equation by
$e^{2y_1}e^{2y_2}$.  Then with
$A_{m_1m_2}=a_{m_1m_2}(e^{-y_1},e^{-y_2})$ and
$C_{m_1m_2}=e^{-2y_1}e^{-2y_2}c_{m_1m_2}(e^{-y_1},e^{-y_2})$,
(\ref{lap}) becomes
\begin{equation*}
e^{-2y_2}(D_1^2-m_1^2)A_{m_1m_2}
+e^{-2y_1}(D_2^2-m_2^2)A_{m_1m_2}=C_{m_1m_2}
\end{equation*}
on the first quadrant in $\mathbb{R}^2$, where $D_j$ stands for
the differential operator $\frac{\partial}{\partial y_j}$, and the
boundary conditions are
\begin{align*}
 A_{m_1m_2}(0,y_2) &= 0, \\
A_{m_1m_2}(y_1,0) &=0.
\end{align*}

We extend $A_{m_1m_2}$ and $C_{m_1m_2}$ by odd reflections in the
variables $y_1$ and $y_2$, labelling the extended functions
$\tilde{A}_{m_1m_2}$ and $\tilde{C}_{m_1m_2}$, respectively, and
we look to  solve
\begin{equation*}
e^{-2|y_2|}(D_1^2-m_1^2)\tilde{A}_{m_1m_2}
+e^{-2|y_1|}(D_2^2-m_2^2)\tilde{A}_{m_1m_2}=\tilde{C}_{m_1m_2}.
\end{equation*}
Let $\chi$ be a smooth compactly supported cutoff function in
$\mathbb{R}^2$, symmetric about the origin, such that $\chi\equiv
1$ in a neighborhood of the origin. Then $\chi\tilde{A}_{m_1m_2}$
satisfies
\begin{equation}
\label{cutoff} e^{-2|y_2|}(D_1^2-m_1^2)\chi\tilde{A}_{m_1m_2}
+e^{-2|y_1|}(D_2^2-m_2^2)\chi\tilde{A}_{m_1m_2}=h,
\end{equation}
where $h$ is a compactly supported, odd function of $y_1$ and
$y_2$, which, when restricted to the first quadrant is
$C^{\infty}$ up to the boundary, and, in a neighborhood of the
origin, is equivalent to $\tilde{C}_{m_1m_2}$.

Upon taking Fourier transforms of (\ref{cutoff}) we obtain
\begin{equation}
\label{orig}
\left((\eta_1^2+m_1^2)e^{-2|D_{\eta_2}|}+(\eta_2^2+m_2^2)e^{-2|D_{\eta_1}|}
\right)\left(\chi \tilde{A}_{m_1m_2}\right)^{\hat{}}=\hat{h},
\end{equation}
where $|D_{\eta_j}|$ is the positive square root of
$-\frac{\partial^2}{\partial \eta_j^2}$ for $j=1,2$.  We intend to
invert the operator
\begin{equation*}
(\eta_1^2+m_1^2)e^{-2|D_{\eta_2}|}+(\eta_2^2+m_2^2)e^{-2|D_{\eta_1}|}.
\end{equation*}
\begin{multline}
\label{K}
\left((\eta_1^2+m_1^2)e^{-2|D_{\eta_2}|}+(\eta_2^2+m_2^2)
e^{-2|D_{\eta_1}|}\right)\left(\chi
\tilde{A}_{m_1m_2}\right)^{\hat{}}=
\\ (\eta^2+m^2)\left(\chi
\tilde{A}_{m_1m_2}\right)^{\hat{}}+\\
\left((\eta_1^2+m_1^2)\left((e^{-2|y_2|}-1)\chi
\tilde{A}_{m_1m_2}\right)^{\hat{}}
+(\eta_2^2+m_2^2)\left((e^{-2|y_1|}-1)\chi
\tilde{A}_{m_1m_2}\right)^{\hat{}}\right)=\\
(\eta^2+m^2)(I-K)\left(\chi \tilde{A}_{m_1m_2}\right)^{\hat{}},
\end{multline}
where $\eta^2=\eta_1^2+\eta_2^2$ and $m^2=m_1^2+m_2^2$, $I$ is the
identity operator, and $K$ is the operator defined by
\begin{equation*}
K\hat{\phi}=
\frac{\eta_1^2+m_1^2}{\eta^2+m^2}\left((1-e^{-2|y_2|})\phi\right)^{\hat{}}
+\frac{\eta_2^2+m_2^2}{\eta^2+m^2}\left((1-e^{-2|y_1|})\phi\right)^{\hat{}}
\end{equation*}
for $\phi\in L^2_0(\mathbb{R}^2)$.

Now let $\chi_1=\chi$, and define cutoff functions, $\chi_j$,
which are symmetric about the origin, for $j=1,2\ldots$, such that
$\chi_j=1$ on $\mbox{supp}\chi_{j-1}$. Also, define $T_0=I$ and
$T_j\phi=\big(K\left(\chi_{j}
T_{j-1}\phi\right)^{\hat{}}\big)^{\check{}}$ for $\phi\in L^2_0$
for $j=1,2, \ldots$. We may assume, after restricting the supports
of the $\chi_j$ if necessary, that the following relations hold
\begin{align*}
&\|T_j\phi\|_2<\|\phi\|_2 \quad \forall \phi\in L^2\mbox{ and }\forall j\in\mathbb{N};\\
&\|T_jA_{m_1m_2}\|_2\rightarrow 0 \quad \mbox{as } j\rightarrow
\infty.
\end{align*}
From (\ref{orig}) and (\ref{K}) we have
\begin{equation}
\label{star}
 (I-K)\left(\chi \tilde{A}_{m_1m_2}\right)^{\hat{}}=\hat{\Phi},
\end{equation}
where $\hat{\Phi}=\frac{\hat{h}}{\eta^2+m^2}$, and from
(\ref{star}) we obtain,
\begin{equation}
\label{telescope}
 \chi_{n+1}T_{n}
A_{m_1m_2}-\chi_{n+2}T_{n+1}A_{m_1m_2}=\chi_{n+1}T_{n}\Phi+s_n,
\end{equation}
where
$s_n=(\chi_{n+1}-\chi_{n+2})T_{n+1}A_{m_1m_2}+\chi_{n+1}\left(K\hat{s}_{n-1}\right)^{\check{}}$
and $s_0=(\chi_2-\chi_1)\Phi$.  Equation \ref{telescope} gives
terms of a telescoping series which converges in $L^2$ since
 $\|\chi_{n+2}T_{n+1}A_{m_1m_2}\|_2\rightarrow 0$ as
$n\rightarrow\infty$.  For any $\epsilon>0$ we may also choose the
$\chi_j$ so that
$\|\chi_{n+1}-\chi_n\|_2<\frac{\epsilon}{2^{n+1}\|A_{m_1m_2}\|_2}$
for $n\geq2$ and
$\|\chi_{2}-\chi_1\|_2<\frac{\epsilon}{2\|\Phi\|_2}$ which implies
$\|\sum_{n=0}^{\infty}s_n\|_2<\epsilon$. Hence, we conclude
\begin{equation}
\label{asexpn} \chi
\tilde{A}_{m_1m_2}\overset{L^2}{=}\sum_{n=0}^{\infty}\chi_{n+1}T_n\Phi.
\end{equation}
\remark{To proceed formally, we may take (\ref{asexpn}) as a
starting point, using (\ref{asexpn}) to define a function
$a_{m_1m_2}(r_1,r_2)$ from the transformations above.  Then it is
easy to show, working backwards, that $v_1$, as defined in
(\ref{vexpn}), gives rise to a function $u_1$ which solves
(\ref{laplace}), (\ref{dirb1}), and (\ref{neub1}).  In fact, using
Lemma \ref{genbi}, we can show the boundary conditions are
satisfied in the classical sense.}

\section{Behavior at the distinguished boundary}
Here we find the singular functions which are in the expansion,
(\ref{asexpn}).  We show, in particular,
\begin{prop}
$\forall N\in\mathbb{N}$, $\exists$ polynomials of degree $N$,
$A_N$, $B_N$, and $C_N$, such that, near the origin, modulo terms
which are in $C^N(\overline{\mathbb{R}_+\times\mathbb{R}_+})$,
\begin{equation*}
 A_{m_1m_2}=A_N\log(y_1^2+y_2^2)+B_N+C_N\arctan\frac{y_1}{y_2}.
\end{equation*}
\end{prop}
In the proof of the proposition we shall make use of functions
constructed in \cite{Eh}.  Let
\begin{equation*}
\Phi_1(y_1,y_2)=-\frac{i}{2}\log(y_1^2+y_2^2)
\end{equation*}
and define $\Phi_{l+1}$ to be the unique solution of the form
\begin{equation*}
 p_1(y_1,y_2) \log(y_1^2+y_2^2) +p_2(y_1,y_2),
\end{equation*}
where $p_1$ and $p_2$ are homogeneous polynomials of degree $2l-2$
in $y_1$ and $y_2$ such that $p_2(y_1,0)=0$, to the equation
\begin{equation*}
\frac{\partial \Phi_{l+1}}{\partial y_2}=\frac{1}{2l}y_2\Phi_l
\end{equation*}
for $l\geq 1$.  Then with $\Phi_l$ defined for $l\geq 1$, define
$(\Phi_l)_0= \Phi_l$ for $y_2\geq 0$, and, for $j \geq 1$,
$(\Phi_l)_j$ to be the unique solution of the form
\begin{equation*}
 p_1\log (y_1^2+y_2^2) + p_2 +p_3\arctan
\left(\frac{y_1}{y_2}\right)
\end{equation*}
 on the half-plane
$\{(y_1,y_2):y_2\geq 0\}$, where $p_1$, $p_2$, and $p_3$ are
polynomials in $y_1$ and $y_2$ such that $p_2(0,y_2)=0$, to the
equation
\begin{equation*}
\frac{\partial(\Phi_l)_j}{\partial y_1}=(\Phi_l)_{j-1}.
\end{equation*}
Also, define recursively for $k\geq 1$, on $y_2\geq 0$,
\begin{equation*}
(\Phi_l)_{jk}=\int_0^{y_2}\cdots\int_0^{t_2}\int_0^{t_1}
(\Phi_l)_j(y_1,t)dt dt_1\cdots dt_{k-1}.
\end{equation*}
\begin{proof}[Proof of the proposition.]
We shall prove that with $T_n$ defined as above,
 $\forall$ $N\in\mathbb{N}$, and $\forall$ $n\geq 0$, on
 $\mathbb{R}_+\times\mathbb{R}_+$, in a neighborhood of $(0,0)$,
\begin{equation}
\label{form} T_n\Phi=\sum_{a+b+2l-2+j+k=2\atop l,j,k\geq 1}^{N}
c_{abljk}y_1^ay_2^b(\Phi_l)_{jk}+s,
\end{equation}
where $c_{abljk}$ depend on $\theta_1$, $\theta_2$, $m_1$, and
$m_2$, and $s$ is used to denote the Fourier transform of any
function which, when restricted to
$\mathbb{R}_+\times\mathbb{R}_+$, is in
$C^{N}(\overline{\mathbb{R}_+\times\mathbb{R}_+})$ (plus terms
which may be singular either along all of $y_1=0$ or along all of
$y_2=0$).
 The proof is by induction.  (\ref{form})
holds true when $n=0$, as shown in \cite{Eh}.  We use the Taylor
expansion with remainder formula,
\begin{equation}
\label{Taylor}
1-e^{-2y_k}=2y_k-\frac{(2y_k)^2}{2!}+\cdots+(-1)^{N}\frac{(2y_k)^{N}}{(N+1)!}
+\frac{(-2)^{N+1}}{(N+1)!}\int_0^{y_k}(y_k-t)^{N+1}e^{-2t}dt,
\end{equation}
for $k=1,2$, in the integrands of the formula
\begin{multline}
\label{nfrom1}
\widehat{T_n\Phi}=\\(-2i)^2\frac{\eta_1^2+m_1^2}{\eta^2+m^2}\int_0^\infty
\int_0^\infty(1-e^{-2y_2})\chi_{n-1}
(T_{n-1}\Phi)\sin(\eta_1y_1)\sin(\eta_2y_2)dy_1dy_2\\
+(-2i)^2\frac{\eta_2^2+m_2^2}{\eta^2+m^2}\int_0^\infty
\int_0^\infty(1-e^{-2y_1})\chi_{n-1}
(T_{n-1}\Phi)\sin(\eta_1y_1)\sin(\eta_2y_2)dy_1dy_2.
\end{multline}
 Now for $y_2 \geq 0$,
\begin{equation}
\label{formofphi} (\Phi_l)_{jk}= p_1\log (y_1^2+y_2^2) + p_2
+p_3\arctan \left(\frac{y_1}{y_2}\right) +p_4\log |y_1|,
\end{equation}
where the $p_m$ are homogeneous polynomials of degree $(2l-2)+j+k$
in $y_1$ and $y_2$ for $m=1,2,3,4$, and we shall also denote by
$(\Phi_l)_{jk}$ its extension to
$\mathbb{R}^2\setminus\{y_1=0,y_2=0\}$, where we use the branch
from 0 to $-\infty$ to extend the $\arctan$ function.

 We show, writing $r_{Nk}$
for the remainder term in (\ref{Taylor}),
\begin{equation}
\label{2bcn} \frac{\eta_i^2+m_i^2}{\eta^2+m^2}\int_0^\infty
\int_0^\infty r_{Nk}\chi_{n-1}
(T_{n-1}\Phi)\sin(\eta_1y_1)\sin(\eta_2y_2)dy_1dy_2,
\end{equation}
for $i=1,2$, is the Fourier transform of a function which may be
included in a function $s$.  We now use the induction hypothesis
so that we may utilize the properties of the particular functions
comprising $T_{n-1} \Phi$.  $r_{Nk}(y_1^ay_2^b(\Phi_l)_{jk})$
vanishes to $(N+2)nd$ order along $y_k$ hence its odd reflection
about the $y_k-$axis will still be $C^{N+1}$ on the appropriate
half-plane. Then the regularity of the operator
$D_1^2+D_2^2-m_1^2-m_2^2$ shows
\begin{equation*}
\frac{\eta_i^2+m_i^2}{\eta^2+m^2}\int_0^\infty \int_0^\infty
r_{Nk}\chi_{n-1}
(y_1^ay_2^b(\Phi_l)_{jk})\sin(\eta_1y_1)\sin(\eta_2y_2)dy_1dy_2
\end{equation*}
is in $C^{N}(\overline{\mathbb{R_+}\times\mathbb{R_+}})$. Again,
using the regularity of $D_1^2+D_2^2-m_1^2-m_2^2$, when the
remaining terms of $T_{n-1}\Phi$ are considered in the integral in
(\ref{2bcn}), we can show that (\ref{2bcn}) may be included in a
function $s$.

After using the induction hypothesis in (\ref{nfrom1}), we
consider
\begin{equation*}
\widehat{\Psi}=\frac{\eta_i^2+m_i^2}{\eta^2+m^2}\int_0^{\infty}\int_0^{\infty}\chi_{n-1}
y_1^py_2^q(\Phi_l)_{jk}\sin(\eta_1y_1)\sin(\eta_1y_1)dy_1dy_2.
\end{equation*}
Instead of looking at the odd function, $\Psi$, of both variables,
we extend $\Psi|_{\mathbb{R}_+\times\mathbb{R}_+}$, denoting the
extended function $\widetilde{\Psi}$, in such a way that
\begin{equation*}
\widehat{\widetilde{\Psi}}=\frac{\eta_i^2+m_i^2}{\eta^2+m^2}
\left(\chi_{n-1}y_1^py_2^q(\Phi_l)_{jk}\right)^{\hat{}}.
\end{equation*}

Then using the relations
\begin{align*}
\frac{\partial}{\partial y_1}(\Phi_l)_{jk}&=(\Phi_l)_{(j-1)k},\\
\frac{\partial}{\partial y_2}(\Phi_l)_{jk}&=(\Phi_l)_{j(k-1)},\\
\frac{\partial}{\partial y_1}\Phi_l&=y_1\Phi_{l-1},\\
\frac{\partial}{\partial y_2}\Phi_l&=y_2\Phi_{l-1},\\
\widehat{\chi_{n-1}(\Phi_l)}_{jk}&=\frac{1}{\eta_1^j}\frac{1}{\eta_2^k}\frac{1}
{(\eta_1^2+\eta_2^2)^l}+s,
\end{align*}
where
\begin{equation*}
\Phi_0=\frac{1}{y_1^2+y_2^2},
\end{equation*}
 we may write $\widehat{\widetilde{\Psi}}$ as a sum of terms
of the form
\begin{equation*}
\left(\varphi
y_1^{\alpha}y_2^{\beta}(\Phi_l)_{jk}\right)^{\hat{}}+s,
\end{equation*}
where $\varphi\in C^{\infty}_0$ is equivalent to 1 in a
neighborhood of the origin, $\alpha$ and $\beta$ are non-negative
integers and $l$, $j$, and $k$ are positive integers.

Once (\ref{form}) is proved, another induction argument shows
$$T_n\Phi\big|_{\mathbb{R}_+\times\mathbb{R}_+}\in
C^{n}(\overline{\mathbb{R}_+\times\mathbb{R}_+})$$ (modulo a
function, $s$),
 and thus we
may prove the proposition by looking at only the first $N$ terms
in (\ref{asexpn}), using Lemma \ref{genbi} to argue the vanishing
of singular terms arising from (\ref{form}) or (\ref{formofphi}).
\end{proof}

%The proof will be finished if we can show terms of the form
%\begin{equation*}
%\frac{Q_l}{(\eta_1^2+\eta_2^2)^l}\hat{\varphi}_j
%\end{equation*}
%are transforms of functions in $C^N(\mathbb{R}^2)$ whenever
% $l\geq N$.  Let $\chi_\eta(\eta_1,\eta_2)$ be a smooth function with the properties
%$\chi_\eta=1$ when $\eta_1^2+\eta_2^2<a$ and $\chi_\eta=0$ when
%$\eta_1^2+\eta_2^2>b$ for some $0<a<b$.  Also, let
%$\chi_\eta'=1-\chi_\eta$.  Then
%\begin{equation*}
%\frac{Q_l}{(\eta_1^2+\eta_2^2)^l}\hat{\varphi}_j=
%\chi_\eta\frac{Q_l}{(\eta_1^2+\eta_2^2)^l}\hat{\varphi}_j+
%\chi_\eta'\frac{Q_l}{(\eta_1^2+\eta_2^2)^l}\hat{\varphi}_j.
%\end{equation*}
%With
%\begin{equation*}
%\phi_1=\chi_\eta\frac{Q_l}{(\eta_1^2+\eta_2^2)^l}\hat{\varphi}_j,
%\end{equation*}
%we have
%\begin{equation*}
%\eta_1^{\alpha}\eta_2^{\beta}(\eta_1^2+\eta_2^2)^l\phi_1\in L^2
%\end{equation*}
%$\forall$ $\alpha$, $\beta\geq 0$, and since
%$(\eta_1^2+\eta_2^2)^l$ is the symbol of an elliptic operator,
%$\phi_1\in C^{\infty}(\mathbb{R}^2)$.  That
%$\chi_\eta'\frac{Q_l}{(\eta_1^2+\eta_2^2)^l}\hat{\varphi}_j$ is in
%$C^N(\mathbb{R}^2)$ follows from
%\begin{equation*}
%\eta_1^{\alpha}\eta_2^{\beta}
%\chi_\eta'\frac{Q_l}{(\eta_1^2+\eta_2^2)^l}\hat{\varphi}_j\in
%L^2(\mathbb{R}^2)
%\end{equation*}
%$\forall$ $0\leq\alpha$, $\beta\leq N$.

After using the decay of $c_{m_1m_2}$ with respect to $m_1$ and
$m_2$ to sum over $m_1$ and $m_2$, and then transforming back to
the variables $z_1$ and $z_2$, we deduce that, $\forall n\in
\mathbb{N}$ $v_1$ may be written
\begin{equation*}
v_1=a_n\log\big(\log r_1)^2+(\log r_2)^2\big) + b_n
+c_n\arctan\left(\frac{\log r_1}{\log r_2}\right),
\end{equation*}
where $a_n$, $b_n$, and $c_n$ are polynomials of degree $n$ in
$\log r_1$ and $\log r_2$ whose coefficients are smooth functions
of $\theta_1$ and $\theta_2$.

We now obtain the singularities of $u_1$ from those of $v_1$. If
\begin{equation*}
u_1=2\sum_{m_1,m_2=-\infty}^{\infty}
b_{m_1m_2}(r_1,r_2)e^{im_1\theta_1} e^{im_2\theta_2},
\end{equation*}
then $u_1$ and $v_1$ are related by
\begin{equation}
\label{a2b}
 \frac{\partial}{\partial
r_2}b_{m_1m_2}-m_2\frac{b_{m_1m_2}}{r_2}=a_{m_1(m_2+1)}.
\end{equation}
We assume without loss of generality that, $\forall$ $m_1$,
$m_2\geq 0$, $b_{m_1m_2}$ and $a_{m_1m_2}$ are supported in some
neighborhood of $r_1=r_2=1$. The solution to (\ref{a2b}) is given
by
\begin{equation*}
b_{m_1m_2}=r_2^{m_2}\int_0^{r_2}t^{-m_2}a_{m_1(m_2+1)}(r_1,t)dt.
\end{equation*}
We make the substitution $u=-\log t$ in the above integral to get
\begin{equation}
\label{int2log}
 b_{m_1m_2}=r_2^{m_2}\int_{-\log
r_2}^{\infty}e^{u(m_2+1)}a_{m_1(m_2+1)}(r_1,e^{-u})du.
\end{equation}
The integral in (\ref{int2log}) was considered in \cite{Eh}
 and gives, after summing over $m_1$ and $m_2$, and using a theorem of Borel,
 with similar
results on the form of $u_2$, Theorem \ref{main}.

We note that there are $f\in C^{\infty}(\overline{\Omega})$, for
example those $f$ whose components, $f_1$ and $f_2$, are
equivalently equal to 1 in a neighborhood of $\partial
\mathbb{D}_1\times\partial\mathbb{D}_2$, which make Theorem
\ref{main} non-trivial, i.e. $\alpha_j$ and $\gamma_j$ are not
necessarily 0.

We may also determine a sufficient condition under which the
solution exhibits any desired degree of regularity up to the
boundary of the bi-disc.
\begin{prop}
If
\begin{equation}
\label{suff}
 \left.\frac{\partial^{2j}}{\partial
r_1^{2j}}\frac{\partial^{2k}}{\partial r_2^{2k}}
\left(\frac{\partial f_1}{\partial
\bar{z}_2}\right)\right|_{r_1=r_2=0}=0
\end{equation}
$\forall$ $j,k\geq 0$ such that $j+k\leq n+2$, then $u_1\in
C^n(\overline{\Omega})$.
\end{prop}
\begin{proof}
If (\ref{suff}) holds, then $\forall$ $m_1, m_2$
\begin{equation*}
\left.\frac{\partial^{2j}}{\partial
y_1^{2j}}\frac{\partial^{2k}}{\partial y_2^{2k}}
C_{m_1m_2}\right|_{y_1=y_2=0}=0,
\end{equation*}
$\forall$ $j,k\geq 0$ such that $j+k\leq n+2$, which implies
$A_{m_1m_2}\in C^n(\overline{\mathbb{R}_+\times\mathbb{R}_+})$
(see \cite{Eh}), $A_{m_1m_2}$ and $C_{m_1m_2}$ defined as above,
and thus $v_1=\frac{\partial u_1}{\partial \bar{z}_2}\in
C^n(\overline{\Omega})$.  Then, we can see $u_1$ is in
$C^n(\overline{\Omega})$ by considering integrals as in
(\ref{int2log}), where now the integrands are in
$C^n(\overline{\Omega})$.
\end{proof}

\end{document}